\documentclass[11pt]{article}

\usepackage{fullpage}
\usepackage{graphicx}
\usepackage{amssymb}
\usepackage{authblk}

\begin{document}
\title{Computing the smallest eigenvalue of large ill-conditioned\\
       Hankel matrices}
\author[1]{Niall Emmart}
\author[2]{Charles C. Weems}
\author[3]{Yang Chen}
\affil[1,2]{Computer Science Department\authorcr University of Massachusetts\authorcr Amherst, MA 01002, USA.\authorcr E-Mail: {\tt nemmart@gmail.com} and {\tt weems@cs.umass.edu}\authorcr\ \ }
\affil[3]{Department of Mathematics\authorcr Imperial College\authorcr 180 Queen's Gate\authorcr London SW7 2BZ, UK \authorcr E-Mail: {\tt y.chen@ma.ic.ac.uk}}
\date{}
\maketitle

\begin{abstract}
This paper presents a parallel algorithm for finding the smallest eigenvalue of a particular form of ill-conditioned Hankel
matrix, which requires the use of extremely high precision arithmetic. Surprisingly, we find that commonly-used approaches
that are designed for high efficiency are actually less efficient than a direct approach under these conditions. We then 
develop a parallel implementation of the algorithm that takes into account the unusually high cost of individual arithmetic 
operations. Our approach combines message passing and shared memory, achieving near-perfect scalability and high tolerance 
for network latency. We are thus able to find solutions for much larger matrices than has been previously possible, with the
potential for extending this work to systems with greater levels of parallelism. 
\end{abstract}

\footnotetext[1]{This work is in supported in part by the National Science Foundation under
NSF grant \#CNS-0619337. Any opinions, findings conclusions or recommendations 
expressed here are the author(s) and do not necessarily reflect those of the sponsors.}

\section*{Introduction}

In this paper, we study the problem parallelizing the computation of the smallest eigenvalue of a family
of extremely ill-conditioned large Hankel matrices using a Beowulf cluster.  The matrices are $N$ by $N$
and are dense and symmetric, given by the formula:

\[
M_{i,j} = \frac{1}{\beta} \Gamma \Big( \frac{1+i+j}{\beta} \Big) \mbox{\ \ \ \ \ (i, j = 0, 1, 2, ..., N-1)}
\]

Where $\beta$ and $N$ are the parameters that determine $M$ and therefore its eigenvalues.

The gamma function, $\Gamma(x)$, is a continuous extension of the factorial function, namely, $\Gamma(n)=(n-1)!$ 
when $n \in \mathbb{N}$.  The following are two example matrixes for $N=4$, $\beta=1$ and $N=4$,
$\beta=\frac{7}{4}$:

\[
  \left[ \begin{array}{cccc}
    0! & 1! & 2! & 3! \\
    1! & 2! & 3! & 4! \\
    2! & 3! & 4! & 5! \\
    3! & 4! & 5! & 6! \\
  \end{array}
  \right]
  \mbox{\ \ \ and\ \ \ }
  \left[\begin{array}{cccc}
    \frac{4}{7}\Gamma(\frac{4}{7}) & \frac{4}{7}\Gamma(\frac{8}{7}) & \frac{4}{7}\Gamma(\frac{12}{7}) & \frac{4}{7}\Gamma(\frac{16}{7}) \\
    \frac{4}{7}\Gamma(\frac{8}{7}) & \frac{4}{7}\Gamma(\frac{12}{7}) & \frac{4}{7}\Gamma(\frac{16}{7}) & \frac{4}{7}\Gamma(\frac{20}{7}) \\
    \frac{4}{7}\Gamma(\frac{12}{7}) & \frac{4}{7}\Gamma(\frac{16}{7}) & \frac{4}{7}\Gamma(\frac{20}{7}) & \frac{4}{7}\Gamma(\frac{24}{7}) \\
    \frac{4}{7}\Gamma(\frac{16}{7}) & \frac{4}{7}\Gamma(\frac{20}{7}) & \frac{4}{7}\Gamma(\frac{24}{7}) & \frac{4}{7}\Gamma(\frac{28}{7}) \\
  \end{array}\right]
\]

\noindent
The matrices are Hankel moment matrices with the weight function $w(x)$=exp($-x^\beta$).

This work is a follow on to previous work by Yang Chen and Nigel Lawrence \cite{chenlawrence}, who investigated the
asymptotic behavior of the smallest eigenvalue of $M$ as $N\rightarrow\infty$.   In the numeric portion of
their paper, they were able to compute the smallest eigenvalue for matrices up to size 300 by 300.  Their
research showed the problem required both extreme precision and a large amount of computation.  Parallel
computing capability has increased substantially since their initial work in 1999, but even today we estimate
it would take over 250 hours to solve a medium size 1000 by 1000 instance on a uniprocessor and much longer
for large matrices.  Clearly an efficient parallel solution is needed.

Unfortunately, because $M$ is so ill-conditioned, it cannot be solved using standard double precision
arithmetic. Neither LAPACK nor the parallel ScaLAPACK offer sufficient precision.  Instances beyond 25 by
25 can only be solved using an extended precision arithmetic package.  For large matrices, the
intermediate computations must be carried out with over 20,000 bits of precision.

The enormous precision required for the intermediate computations presents opportunities and challenges.
It shifts the balance of computation to communication by requiring much more computation for each scalar 
operation on the matrix elements.  It also means that the
amount of main memory needed to store the matrix is massive. The extended precision computations also 
affect cache locality. Therefore the array elements must be more finely partitioned
among the nodes than normal, with respect to traditional heuristics for the blocking factors of large matrix computations. 

This paper explores various algorithms to solve the problem with the goal of maximizing the delivered
performance from a Beowulf cluster.

\section*{Hankel matrices, the moment problem and smallest eigenvalues}

Hankel matrices occur naturally in moment problems - for a given sequence of moments 
(mean, variance, skewness, kurtosis, etc), determine a probably distribution/measure that 
gives rise to the sequence of moments.  See, for example, the monographs by Akhiezer \cite{akh} and
by Krein \cite{krein} for in-depth studies.  Moment problems are encountered in many fields from 
statistics to quantum physics to hydrology.  They are classified according to 
the support of the distribution/measure (the set of points where the distribition/measure 
is non-zero).  If the support is a closed interval (the Hausdorff moment problem) there
will always be a unique solution.  If the support is the whole number line (the Hamburger moment
problem), then one can encounter a situation where the problem is said to be indeterminate,
that is there are infinitely many probability distributions/measures with the same sequence
of moments.

Hankel matrices also appear in random matrix theory, see Mehta \cite{mehta}, a currently active
research area that encompasses pure and applied mathematics and has applications in
wireless communications and multi-variate statistics.  See for example the monograph by
Kerov \cite{kerov} which studies the connections between representation theory, moment problems
and random matrices.  

It was found in the classical papers by Szeg\"{o} \cite{szego} and by Widom and Wilf \cite{wandw}
that for those Hankel matrices (of order $n$) generated by a probability density supported on
a closed interval, the corresponding smallest eigenvalue tends to zero as $n$ tends to infinity.

In \cite{bci} Berg, Chen and Ismail showed that in the infinite interval case, 
the moment problem will be indeterminate if and only if the smallest eigenvalue of the Hankel matrix 
tends to a strictly positive number as $n$ tends to infinity.  This is a new criteria for the 
determinancy of the Hamburger moment problem and hence our motivation for fast algorithms
to compute the smallest eigenvalues of large Hankel matrices.

\section*{Properties of M}

The algorithms we will explore in this paper rely on three properties of $M$:

\begin{itemize}
\item $M$ is real and symmetric, therefore, all the eigenvalues are real.
\item The eigenvalues of $M$ are the zeros of the characteristic polynomial $P(x)=\det(M-xI)$.
\item $M$ is the Hankel moment matrix for the weight function $w(x)$=exp($-x^\beta$), therefore,
there are $N$ distinct positive eigenvalues of $M$.
\end{itemize}

The last point is very important for convergence of the Secant algorithm and can be proven as 
follows.  Since $w(x)$=exp($-x^\beta$), we know that for any real polynomial $Q(x)$ of degree
$N$, the moment integrals ($i=0,1,2,...$) will all exist:

$$
\mu_{i}:=\int_{0}^{\infty} x^{i}w(x)dx 
$$
  
further, since 
 
$$
 \int_{0}^{\infty} [Q(x)]^2 w(x)dx>0
$$

the associated symmetric Hankel matrix ${\cal H}_N:=(\mu_{i+j})_{i,j=0}^{N},$ will be
positive definite and the eigenvalues will be distinct for any fixed $N$.  See {\it
Matrix Analysis,} Horn and Johnson, Cambridge University Press, 1985 \cite{hornjohnson}.

\section*{M is extremely ill-conditioned}

The condition number of a matrix is a good indication of how much precision is needed to
compute the smallest eigenvalues.  The larger the condition number the greater the required precision.

The condition number is defined to be:
\[
cond(M) = \frac{\lambda_N[M]}{\lambda_1[M]}  \mbox{\ \ \ \  $\lambda_N$ is the largest eigenvalue and $\lambda_1$ is the smallest}
\]

Condition numbers can be estimated using the computation of $\lambda_1$ and by bounding
$\lambda_N$ as follows.  The Raleigh quotient function, $\rho(u;M)$
ranges over the interval $[\lambda_1, \lambda_N]$ for non-zero $u$, therefore:

\[
  \rho(u;M) = \frac{\langle Mu, u \rangle}{\langle u, u \rangle} \le \lambda_N \mbox{\ \ \ \ \ For all non-zero $u$}
\]

Choosing $u$ to be the column vector $(0, 0, 0, ..., 1)$:

\[
  \rho(u;M) = \frac{\langle Mu, u \rangle}{\langle u, u \rangle} = \frac{1}{\beta}\Gamma \Big(\frac{2N-1}{\beta} \Big) \le \lambda_N
\]

Further, $\lambda_N \le \mbox{trace}(M)$, therefore:

\[
 \frac{1}{\beta}\Gamma \Big(\frac{2N-1}{\beta} \Big) \le \lambda_N \le \mbox{trace}(M)=\frac{1}{\beta}\Gamma \Big(\frac{2N-1}{\beta} \Big) + \sum_{k=0}^{N-2}{\frac{1}{\beta}\Gamma \Big(\frac{2k+1}{\beta} \Big)}
\]

For large $N$, $\frac{1}{\beta}\Gamma \Big(\frac{2N-1}{\beta} \Big)$ is a good estimate for $\lambda_N$ because $\sum_{k=0}^{N-2}{\frac{1}{\beta}\Gamma \Big(\frac{2k+1}{\beta} \Big)}$ is very small compared to $\frac{1}{\beta}\Gamma \Big(\frac{2N-1}{\beta} \Big)$.

The following table of condition numbers has been calculated using the lower bounds for $\lambda_N$ and the results of the $\lambda_1$ computations:

\small\begin{tabular}{|r||l|l|l|l||l|l|l|}
\multicolumn{8}{c}{} \\
\multicolumn{1}{c}{} & \multicolumn{4}{c}{\textbf{LOWER BOUND ON COND(M)}} & \multicolumn{3}{c}{\textbf{EXPONENTIAL FUNCTIONS}} \\
\hline

N        & $\beta$=1/3   & $\beta$=1/2   & $\beta$=1 & $\beta$=7/4 & $2^N$  & $N!$ & $N^N$ \\
\hline
100      & 8.52 e1396  & 7.36 e861   & 9.40 e384   & 1.94 e228  & 1.27 e30  & 9.33 e157  & 1.00 e200  \\
300      & 5.17 e5066  & 4.65 e3167  & 6.38 e1429  & 3.55 e819  & 2.04 e90  & 3.06 e614  & 1.37 e743  \\
500      & 4.56 e9116  & 6.89 e5726  & 3.62 e2597  & 4.11 e1472 & 3.27 e150 & 1.22 e1134 & 3.05 e1349 \\
1000     & 1.85 e20050 & 6.97 e12663 & 7.62 e5780  & 6.80 e3240 & 1.07 e301 & 4.02 e2567 & 1.00 e3000 \\
1500     & 1.11 e31666 & 3.81 e20055 & 8.45 e9187  & 3.32 e5125 & 3.51 e451 & 4.81 e4114 & 1.37 e4764 \\
\hline
\multicolumn{8}{c}{} \\
\end{tabular}\normalsize

As can be seen from this table, in each case the condition number is growing faster than $N^N$.

\section*{Algorithm Selection}

There are several standard techniques for finding eigenvalues, see \cite{parlett}.  For this paper
we implemented four algorithms and evaluated them with a number of criteria:  (a) how much precision
does the algorithm require to meet a desired level of precision in the output;  (b) how fast
is the calculation relative to the other algorithms;  (c) how effectively can the algorithm be parallelized.
We implemented the following algorithms using the GNU Multiple Precision \cite {GMP} library.

\begin{itemize}
\item Lanczos Algorithm as described by Stoer and Bulirsch \cite{stoerbulirsch} p. 401
\item Householder's Algorithm as described by Stoer and Bulirsch \cite{stoerbulirsch} p. 391 (with minor corrections)
\item The Jacobi Method with Rutishauser's enhancements, as described by Parlett \cite{parlett}
      pp. 189-196.
\item A direct approach with an $LDL^T$ determinant algorithm and a Secant root finder, as
      described below (hereafter referred to as the Secant algorithm)
\end{itemize}

The algorithms were tested by varying the number of bits of precision (K) of the inputs and
the computations to achieve a specific level of precision in the result.  In the following
table, {\it Accuracy} is the number of correct significant digits in the result.  {\it Run Time}
is in seconds.  {\it Factor} is the number of times slower this algorithm is than the fastest
for a given $N$.  The algorithm with the most accuracy for the least run time will be the best 
choice for computing the smallest eigenvalue.

\small\begin{tabular}{|l|r|r|r|r|r|r|}
\multicolumn{7}{c}{} \\
\hline
\multicolumn{1}{|c|}{Algorithm} & \multicolumn{1}{c|}{N} & \multicolumn{1}{c|}{$\beta$} & \multicolumn{1}{c|}{K} & \multicolumn{1}{c|}{Accuracy} & \multicolumn{1}{c|}{Run Time} & \multicolumn{1}{c|}{Factor} \\
\hline

Secant      & 100 &    1    & 800    &  60      & 0.96               & 1.0 \\
Householder & 100 &    1    & 2490   &   9      & 3.80               & 3.96 \\
Jacobi      & 100 &    1    & 700    &  58      & 5.91               & 6.15 \\
Lanczos     & 100 &    1    & 206250 &  36      & 361.46             & 376.50 \\

\hline

Secant      & 200 &    1    & 1400   &  90      & 23.67              & 1.0 \\
Householder & 200 &    1    & 5940   &  59      & 134.32             & 5.67 \\
Jacobi      & 200 &    1    & 1500   &  24      & 216.55             & 9.15 \\

\hline

Secant      & 300 &    1    & 1600   &  60      & 144.02             & 1.0 \\
Householder & 300 &    1    & 10000  &  53      & 990.26             & 6.88 \\
Jacobi      & 300 &    1    & 2400   &  8       & 2027.72            & 14.08 \\

\hline

Secant      & 400 &    1    & 2100   & 51       & 603.84            & 1.0 \\
Householder & 400 &    1    & 13250  & 41       & 4036.11           & 6.68 \\
Jacobi      & 400 &    1    & 3375   & 48       & 9888.59           & 16.38 \\
\hline
\multicolumn{7}{c}{} \\
\end{tabular}\normalsize

One unexpected result was the poor performance of the Lanczos algorithm, which failed
for N=100 with anything less than 205000 bits of precision.  There are a number of variants of the basic
Lanczos algorithm - some orthogonalize the new vector after each iteration and some stop and 
restart the process, however we did not explore these because it did not appear that any of them 
would reduce the precision requirements to less than that of the other three algorithms.

From these results we can conclude that the fastest technique to compute the smallest
eigenvalue under these conditions is the Secant algorithm. In the remainder of the paper, we show that it
can be effectively parallelized.

\section*{The Secant Algorithm}

The Secant algorithm can be used to find the smallest root of the characteristic 
polynomial, $P(x)$.  $P(x)$ can be defined as either $\det(xI-M)$ or $\det(M-xI)$. 
We use the latter because it guarantees $P(0)$ will be positive.  The secant
algorithm starts with two initial points $x_1$ and $x_2$ which must be less 
than $\lambda_1$.  We choose $x_1$ to be a small negative value and $x_2$ to
be zero.  The Secant algorithm then computes a sequence of $x_i$'s using the 
following recurrence relation:

\[
  x_{i+1} = x_i - \frac{x_i - x_{i-1}}{P(x_i) - P(x_{i-1})}P(x_i)
\]

The computation is complete when the most significant digits (15 decimal digits 
in our case) of $x_i$ have stabilized.  Since the eigenvalues of $M$ are all unique,
the roots of $P$ are all simple and the Secant algorithm is guaranteed to 
converge rapidly to $\lambda_1$ \cite{burdenfaires}.  Also note, since $P$ has no inflection points
less than $\lambda_1$, the slope of the Secant at $x_i$ will always be less 
than $P'(x)$ when $x_i \le x \le \lambda_1$ and therefore the Secant is
guaranteed not to overshoot $\lambda_1$.  

Finally, there is no need to solve for $P$, instead, we can evaluate 
$P(x)$ directly by computing $\det(M-xI)$.  Thus, for our matrices, the problem of finding
the smallest eigenvalue reduces to that of solving a sequence of determinants.

\section*{Interval Verification}

After the secant method has converged to some $x$, the next step is to prove that $x$ is
in fact an eigenvalue of $M$.  To do this, we do two checks, first we truncate $x$ to 15
significant digits and compute $\det(M-xI)$ which must be greater than 0.  Then we add 1
to the least significant digit of $x$ and compute $\det(M-xI)$ again.  This result must be
less than 0.  To eliminate the possibility that round-off errors have poisoned the result,
we perform these checks using fixed point interval arithmetic.   See, for example, {\it
Interval Analysis,} R.E. Moore \cite{remoore}. In general the verification algorithm requires 
far more precision than the Secant root finder. 

\section*{The $LDL^T$ Determinant Algorithm}

There are several standard techniques to compute the determinant of a matrix.
The fastest general methods all involve factoring the matrix in some way and then
calculating the determinant from the factors.  In our case, we can make
use of the fact that $M$ is symmetric and factor $M-xI$ into a lower
triangular matrix, $L$, a diagonal matrix, $D$, and the transpose, $L^T$,
such that $M-xI=LDL^T$.  Because $x<\lambda_1$, $M-xI$ will be positive
definite and the factorization is essentially the same as a Cholesky
factorization except that it avoids the square root operations.  The
$LDL^T$ factorization has numerical stability comparable to Cholesky
factorization.  Once the matrix is factored, the determinant is simply
the product of the elements on the main diagonal of $D$.

To make the algorithm easy to parallelize, we use the `submatrix' order 
\cite{demmel} for the $LDL^T$ algorithm which applies a column of the matrix to 
all the remaining columns to its right.  Presented below are serial and parallel 
versions of the algorithm.  For clarity they both are shown using double precision 
arithmetic.  In the actual implementations the calculations are done using
the GMP integer package with fixed-point numbers.  The elements of $M$ are 
stored with $K$ bits after the decimal point, which is specified when the 
algorithm is started.  The elements of C, nextC, CDivDiag, and nextCDivDiag
are all stored with K/2 bits after the decimal point.

\subsubsection*{Serial $LDL^T$ Algorithm:}

\small\begin{verbatim}
double LDLTDeterminant(double M[][], double x) {
  int    processingCol, row, col, N=size(M);
  double CDivDiag[N], determinant=1.0;
     
  for(processingCol=0; processingCol<N; processingCol++) {
    M[processingCol][processingCol] -= x;
    determinant = determinant * M[processingCol][processingCol];
    
    // CDivDiag is set to the current column that we're processing divided by
    // the entry on the diagonal.
    for(row=processingCol+1; row<N; row++)
      CDivDiag[row] = M[row][processingCol] / M[processingCol][processingCol];
    for(col=processingCol+1; col<size; col++)
      for(row=col; row<N; row++)
        M[row][col] -= M[col][processingCol] * CDivDiag[row];
  }
  return determinant;
}
\end{verbatim}\normalsize

\subsubsection*{Parallel $LDL^T$ Algorithm:}

The parallel version of the $LDL^T$ algorithm is designed to use both MPI and
OpenMP.  MPI is used to distribute data between nodes and OpenMP is used to spread
the computation across multiple cores.  We chose this solution for two reasons.
First, the inter-core communication overhead of OpenMP is much less than for MPI.  Second,
the majority of the computation time is spent in just a few loops that naturally parallelize
with OpenMP.  The data distribution is done using MPI broadcasts.  These are run in
a separate thread which allows the communication to be overlapped with the computation.

The following pseudo-code has annotations on the right hand side indicating aspects of 
our timing analysis, for example, \{c\}. 
The algorithm tracks the total amount of time spent executing \{c\} steps, likewise for \{net\}
and \{d\}.  $c$ is for the main computation, $net$ is for the time spent waiting 
for the network IO to complete, and $d$ for the remaining divisions required to
complete CDivDiag.  In addition one more timer records the total time spent 
computing the determinants.

\small\begin{verbatim}
void assign(int assignments[]) {
  // This routine assigns the columns to MPI processes.  A column is assigned when 
  // assignments[col]=rank (the MPI rank of the process).  N is the number of columns
  //  S is MPI size.
  for(col=0; col<N; col++) {
    rank = col % (2*S)
    assignments[col] = (rank<S) ? rank : (2*S)-1 - rank;
  }
}

void apply(int processingCol, double C[], double CDivDiag[], double[][] M, 
           int assignments[], int startCol, int endCol) {
  // This routine applies the column C to columns of M from start to end
  for(col=startCol; col <= endCol; col++)
    if(assigments[col]==myRank)
      for(row=col; row<N; row++)
        M[row][col] -= C[processingCol] * CDivDiag[row]
}

void startBackgroundTransmit(Column C, Column CDivDiag) {
  // This routine is run by a background thread.  It uses MPI_Bcast to 
  // transmits the column C to all of the other MPI processes. 
  // While there is time left in the foreground computation, 
  // send the next chunk of the CDivDiag column.
}

void startBackgroundReceive(Column C, Column CDivDiag) {
  // Receive the C column.
  // Receive any values of the CDivDiag column that are sent.  If the 
  // computation is IO bound, no values will be sent.  If it is CPU 
  // bound, all the CDivDiag values will be sent.
}

int waitForBackgroundCommunication() {
  // Wait until the background IO is finished
  // Return the number of CDivDiag entries sent/received
}

double determinant(double M[][], double x) {
   double C[N], CDivDiag[N], nextC[N], nextCDivDiag[N], determinant=1.0;
   int    processingCol, row, col, CDivDiagReceived, assignments[N];

   // Compute the assignments
   assign(assignments);

   // Set M to M-xI   
   for(row=0; row<N; row++)
    M[row][row] -= x;    
   
   // Pull out the first column, and compute CDivDiag
   copyColumn(M, 0, C);
   for(row=1; row<N; row++)
    CDivDiag=C[row] / C[0];
    
   For(processingCol=0; processingCol<N-1; processingCol++) {
     if(assignments[processingCol+1] == myRank) {
       // Apply C to column processingCol+1
       apply(processingCol, C, CDivDiag, M, assignments, processingCol+1,          {c}
             processingCol+1);
             
       // Pull the column out of M.  This column is finished and needs to be
       // sent to the other MPI processes.  Compute nextCDivDiag.
       copyColumn(M, processingCol+1, nextC);
       for(row=processingCol+2; row<N; row++)
         nextCDivDiag[row] = nextC[row] / nextC[processingCol+1];

       startBackgroundTransmit(nextC, nextCDivDiag);       

       // Apply C (processingCol) to the remainder of M
       apply(processingCol, C, CDivDiag, M, assignments, processingCol+2, N-1);    {c}
       waitForBackgroundCommunication();                                           {net}
     }
     else {
       startBackgroundReceive(nextC, nextCDivDiag);

       // Apply C to all remaining columns of M
       apply(processingCol, C, CDivDiag, M, assignments, processingCol+1, N-1);    {c}
       CDivDiagReceived=waitForBackgroundCommunication();                          {net}

       // Compute any nextCDivDiag entries not received
       for(row=processingCol+2+CDivDiagReceived; row<N; row++)
         nextCDivDiag[row] = nextC[row] / nextC[processingCol+1];                  {d}
     }
     swap(C, nextC);
     swap(CDivDiag, nextCDivDiag);
     determinant=determinant * C[processingCol+1];
   }
   return determinant;
}
\end{verbatim}\normalsize

\section*{Numerical results}
The following table shows a sampling of results for the smallest eigenvalue computations. 
It clearly illustrates the wide range of values that the algorithm must handle.

\begin{tabular}{|r|l|l|l|l|}
\multicolumn{5}{c}{} \\
\multicolumn{5}{c}{THE SMALLEST EIGENVALUE OF M} \\
\hline

N        & $\beta$=1/3   & $\beta$=1/2       & $\beta$=1          & $\beta$=7/4        \\
\hline
100      & 3.4720        & 2.7397x$10^{-1}$  & 2.1079x$10^{-15}$  & 1.6976x$10^{-45}$  \\
300      & 3.3984        & 1.5837x$10^{-1}$  & 5.5215x$10^{-28}$  & 1.4844x$10^{-102}$ \\
500      & 3.3763        & 1.2047x$10^{-1}$  & 1.1138x$10^{-36}$  & 6.7121x$10^{-149}$ \\
1000     & 3.3544        & 8.2087x$10^{-2}$  & 1.0892x$10^{-52}$  & 3.6209x$10^{-246}$ \\
1500     & 3.3447        & 6.6295x$10^{-2}$  & 5.4593x$10^{-65}$  & 6.4232x$10^{-330}$ \\
\hline
\multicolumn{5}{c}{} \\
\end{tabular}

\section*{Performance}

The performance testing was done on the University of Massachusetts Swarm cluster, which 
has 60 compute nodes, each with 8 cores (Xeon 5355 processors at 2.66 GHz) and 16 GB of
RAM per node.  The nodes are connected via gigabit ethernet.  The cluster is partitioned 
and our tests were run on 48 nodes, each with 5 cores for a total of 240 cores.

The results show the cumulative time spent running the determinant algorithm.  The results
do not include the startup time, the time to generate the M matrix, nor the time spent
running the secant algorithm.  The total for these tasks is less than 1\% of the time spent
computing determinants.

In the tables and graphs below {\it Total Time} is the total time spent computing determinants.
{\it Computation} is the cumulative time spent executing \{c\} steps in the $LDL^T$ algorithm and
{\it Net + Divs} is the cumulative time in the \{net\} and \{d\} steps.  For these results, we
add the \{net\} and \{d\} times together because they both represent time/computations wasted
due to a lack network bandwidth.

\noindent\small\begin{tabular}{|c|l|r r|r r|r r|}
\multicolumn{8}{r}{} \\
\multicolumn{8}{c}{CONSTANT BETA} \\
\hline
\multicolumn{2}{|r|}{} & \multicolumn{2}{|c|}{$\beta=7/4$, N=1000} & \multicolumn{2}{|c|}{$\beta=7/4$, N=1500} & \multicolumn{2}{|c|}{$\beta=7/4$, N=2000} \\
\hline
 80   & Total Time     & 0:43:57 & (2,638)               & 4:04:31  & (14,672) & 13:57:35  & (50,255) \\
Cores & Computation    & 92.2\%  & (2,432)               & 97.0\%   & (14,231) & 98.4\%    & (49,438) \\
      & Net + Divs     &  7.2\%  &   (190)               &  2.8\%   &    (417) &  1.5\%    &    (773) \\
\hline
160   & Total Time     & \multicolumn{2}{c|}{No results} & 2:04:00  & (7,440)  & 7:05:11   & (25,511) \\
Cores & Computation    & \multicolumn{2}{c|}{available}  & 93.6\%   & (6,965)  & 96.9\%    & (24,727) \\
      & Net + Divs     & \multicolumn{2}{c|}{}           &  6.1\%   &   (452)  &  3.0\%    &    (752) \\
\hline
240   & Total Time     & 0:17:33 & (1054)                & 1:24:53  &  (5,094) & 4:45:25   & (17,125) \\
Cores & Computation    & 77.0\%  &  (811)                & 91.2\%   &  (4,645) & 96.2\%    & (16,466) \\
      & Net + Divs     & 21.5\%  &  (227)                &  8.4\%   &    (426) &  3.7\%    &    (630) \\
\hline
\multicolumn{8}{r}{} \\
\end{tabular}\normalsize

The table shows that as $N$ increases, the percentage of time spent doing useful computation increases.

\begin{figure}[h]
\centering
\includegraphics[viewport=17 531 448 773,clip,scale=0.8]{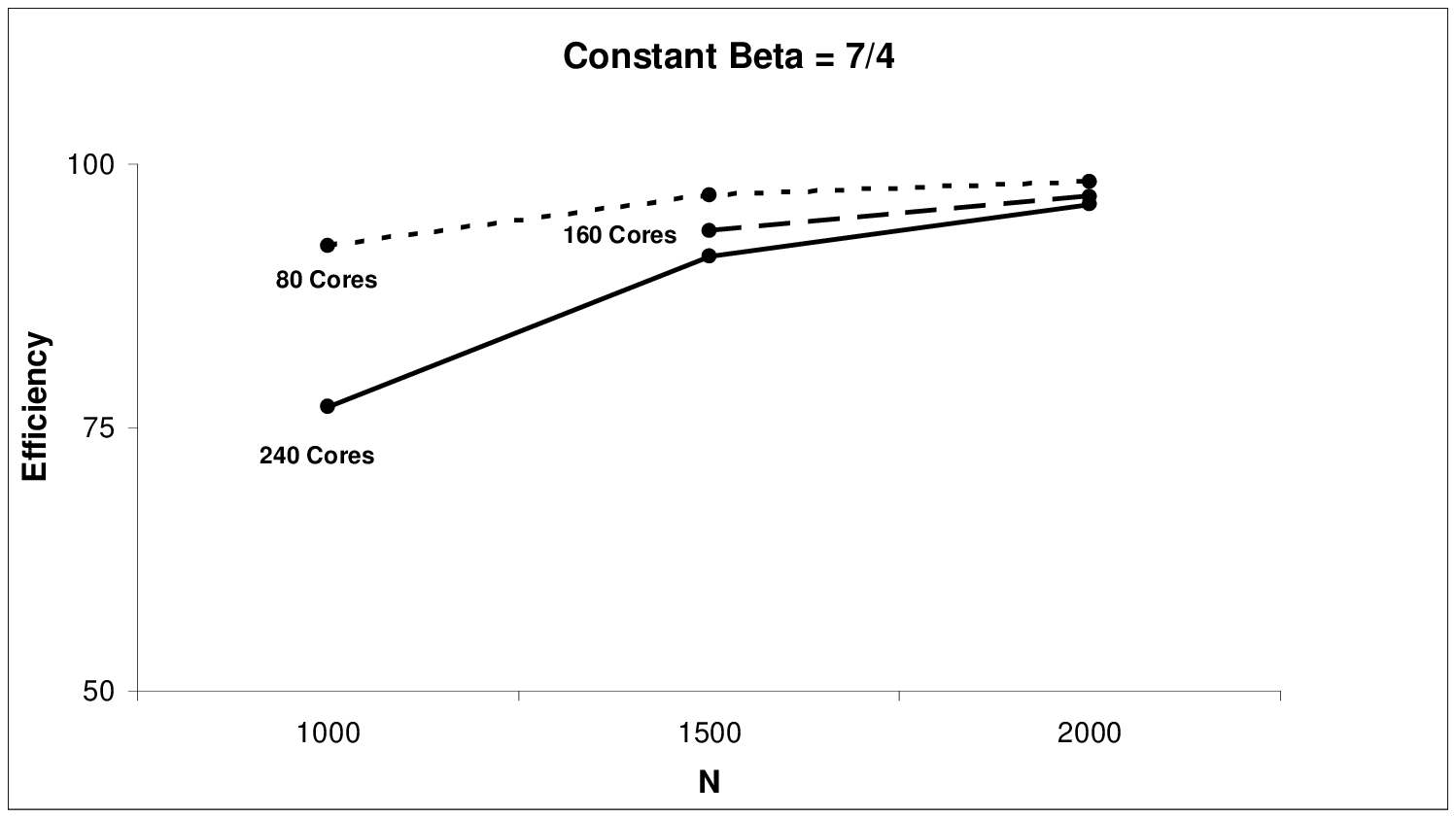}
\caption{Constant Beta}
\end{figure}

\noindent\small\begin{tabular}{|c|l|r r|r r|r r|r r|}
\multicolumn{10}{r}{} \\
\multicolumn{10}{c}{CONSTANT N} \\
\hline
\multicolumn{2}{|r|}{} & \multicolumn{2}{|c|}{$\beta=1/3$, N=1500} & \multicolumn{2}{|c|}{$\beta=1/2$, N=1500} & \multicolumn{2}{|c|}{$\beta=1$, N=1500} & \multicolumn{2}{|c|}{$\beta=7/4$, N=1500} \\
\hline
 80   & Total Time     & 8:29:45 & (30,585) & 6:38:32  & (23,912) & 3:59:40   & (14,380) & 4:04:31  & (14,672) \\
Cores & Computation    & 88.5\%  & (27,075) & 87.3\%   & (20,868) & 85.4\%    & (12,277) & 97.0\%   & (14,231) \\
      & Net + Divs     & 11.4\%  &  (3,483) & 12.6\%   &  (3,013) & 14.4\%    &  (2,076) & 2.8\%    &    (417) \\
\hline
160   & Total Time     & 6:47:22 & (24,442) & 4:10:38  & (15,038) & 3:26:40   & (12,400) & 2:04:00  &  (7,440) \\
Cores & Computation    & 69.7\%  & (17,031) & 68.6\%   & (10,308) & 46.5\%    &  (5,761) & 93.6\%   &  (6,965) \\
      & Net + Divs     & 30.2\%  &  (7,388) & 31.3\%   &  (4,699) & 53.3\%    &  (6,613) &  6.0\%   &    (452) \\
\hline
240   & Total Time     & 4:21:12 & (15,673) & 3:38:29  & (13,109) & 2:23:07   &  (8,588) & 1:24:53  &  (5,094) \\
Cores & Computation    & 54.7\%  &  (8,575) & 49.3\%   &  (6,469) & 43.9\%    &  (3,767) & 91.2\%   &  (4,645) \\
      & Net + Divs     & 45.1\%  &  (7,073) & 50.4\%   &  (6,609) & 55.8\%    &  (4,793) &  8.4\%   &    (426) \\
\hline
\multicolumn{8}{r}{} \\
\end{tabular}\normalsize

\begin{figure}[h]
\centering
\includegraphics[viewport=17 531 448 773,clip,scale=0.8]{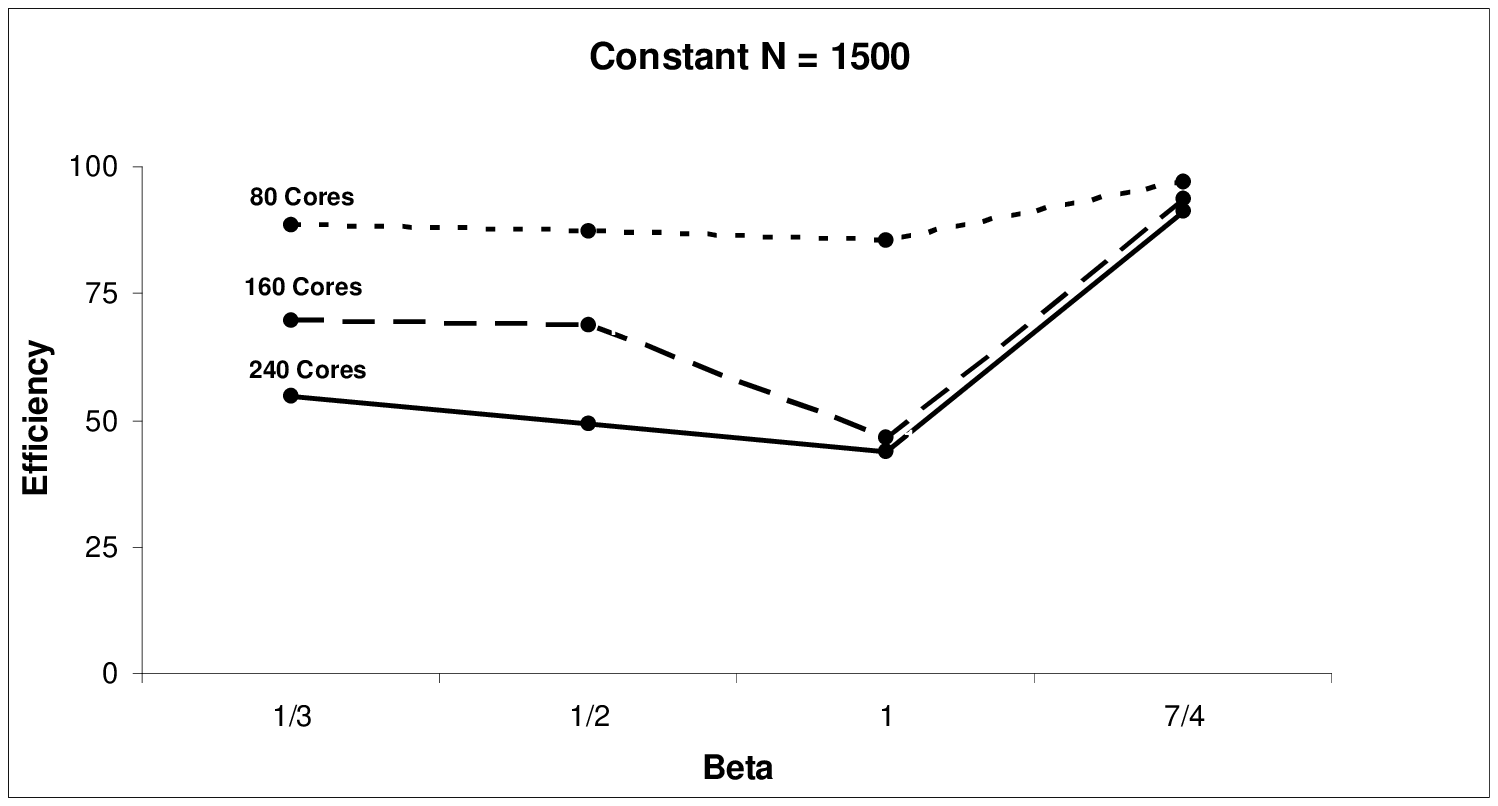}
\caption{Constant N}
\end{figure}

This table shows two trends.  First, as $\beta$ moves away from 1, the computation becomes more efficient.
That is, the percentage of time spent doing computations increases and the percentage of time spend waiting
on network IO and divisions decreases.  Second, if the computation is network bound, increasing the number
of cores makes the computation substantially less efficient.
However, if the computation is CPU bound, then increasing the number of cores makes the computation only
marginally less efficient.  In fact, the amount of time wasted waiting for network IO and divisions is almost
constant.  The percentage of wasted time increases only because the overall time to complete the computation
decreases.  In other words, the algorithm scales almost perfectly when it is CPU bound.

\section*{Communication}

There are a couple of interesting things to note about the algorithm.  First, it works by
preparing a column, starting a background broadcast and then applying the column to the remainder of
the matrix.   As the matrix size grows, or the precision in the numbers grows, it becomes easier and
easier to fully overlap the computation and communication.  This means that, with a large enough problem,
network latency has minimal impact on performance, which is dominated by network throughput.  Given
sufficient network throughput, the algorithm can scale to very large systems.

Second, as $N$ grows, the precision ($K$) required to perform the computation increases.
In the inner loop of the computation two numbers, each with $K/2$ bits are being multiplied. 
Increasing $K$ thus increases the processing time by approximately $O(K^2)$, while the communication time
increases only linearly.  Therefore, as $N$ grows, there is a very powerful effect on the efficiency
and scalability of the algorithm.

\section*{Conclusion}

Large ill-conditioned Hankel matrices present an unusual mix of computation that are not readily 
solved with  traditional approaches. We have explored a space of potential algorithmic solutions to 
determine which approach provides the greatest efficiency given the special precision requirements of 
the problem. Surprisingly, we found that a direct method is more effective than the more sophisticated 
algorithms that are commonly employed. Our parallel implementation of this algorithm takes into 
account the atypically large computation time required for individual operations, and uses both 
shared memory and message passing to optimize performance. The result is a parallel 
implementation that scales nearly perfectly, and that can be made nearly insensitive to network 
latency while taking maximum advantage of network throughput. The algorithms we developed 
have thus proved to be an elegant, efficient, fast and scalable solution to the problem, with 
guaranteed numeric results. As a result, we have been able to considerably extend the known 
set of solutions for these matrices using a modest Beowulf cluster, and have shown that much 
larger instances should be easy to solve on systems with a higher degree of parallelism.


\begin{thebibliography}{1}

  \bibitem{akh}
     N. I. Akhiezer,
     \emph{The classical moment problem and some related questions in analysis.}
     Oliver and Boyd,
     Edinburgh,
     1965.
     
  \bibitem{bci}
     C. Berg, Y. Chen and M.E.H. Ismail
     Small eigenvalues of large Hankel matrices: the indeterminate case.
     \emph{Math. Scand.} vol. 91, 67--81.
     2002.
     
  \bibitem{burdenfaires}
     Richard L. Burden, J.Douglas Faires.
     \emph{Numerical Analysis}, 4th edition.
     PWS-KENT Publishing Company,
     Boston, Massachusetts,
     1989.

  \bibitem{chenlawrence}
     Y. Chen and N. D. Lawrence.
     Small eigenvalues of large Hankel matrices.
     \emph{J. Phys. A: Math. Gen}, vol. 32, 7305--7315.
     1999.

  \bibitem{demmel}
     James W. Demmel, Michael T. Heath, Henk A. van der Vorst.
     Parallel Numerical Linear Algebra.
     ACTA Numerica, 1992.
     
  \bibitem{GMP}
    GNU Open Source Community.
    \emph{The GNU Multiple Precision Arithmetic Library}.
    http://www.gmplib.org/

  \bibitem{hornjohnson}
     Horn and Johnson.
     \emph{Matrix Analysis}.
     Cambridge University Press,
     1985

  \bibitem{kerov}
     S. V. Kerov,
     \emph{Asymptotic representation theory of symmetric group and its application in analysis}.
     AMS, 2003.
     
  \bibitem{krein}
     M. G. Krein and A. A. Nudelman.
     \emph{Markov moment problems and extremal problems}.
     AMS, 1977.
     
  \bibitem{mehta}
     M. L. Mehta.
     \emph{Random Matrices}, 3rd edition.
     Elsevier,
     Singapore,
     2006.
     
  \bibitem{remoore}
     R. E. Moore,
     \emph{Interval Analysis}.
     Prentice-Hall Inc.,
     Englewood Cliffs, New Jersey,
     1966.
   
  \bibitem{parlett}
     Beresford N. Parlett,
     \emph{The Symmetric Eigenvalue Problem}.
     Classics in Applies Mathematics; 20.
     Prentice-Hall Inc.,
     Englewood Cliffs, New Jersey,
     1997

  \bibitem{stoerbulirsch}
     J. Stoer, R. Bulirsch.
     \emph{Introduction to Numerical Analysis}, 3rd edition.
     Texts in Applied Mathematics 12.
     Springer Science + Business Media, LLC,
     New York, New York,
     2002
     
  \bibitem{szego}
     G. Szeg\"{o}.
     On some Hermitian forms associated two geiven curves of the complex plane.
     \emph{Trans. Amer. Math. Soc.} vol. 40, 450--461.
     1936.

   \bibitem{wandw}
     H. Widom and H. Wilf.
     Small eigenvalues of large Hankel matrices.
     \emph{Proc. Amer. Math. Soc.} vol. 17, 338--344.
     1966.

\end{thebibliography}
\end{document}